\documentclass[12pt]{amsart}
\usepackage{graphicx,afterpage,amsfonts,amssymb}

\advance\textheight\topskip
\allowdisplaybreaks

\usepackage{empheq}
\usepackage{upgreek}
\newcommand{\pL}{\lambda}
\newcommand{\pB}{\beta}
\newcommand{\pR}{\rho}
\newcommand{\pT}{\tau}

\newcommand{\leftaction}{\kern1pt{\rightharpoonup}\kern1pt}
\newcommand{\rightaction}{\kern1pt{\leftharpoonup}\kern1pt}

\usepackage{times}
\usepackage[mathscr]{eucal}
\usepackage{mathptmx}

\usepackage{amsthm}

\newcommand{\oCq}{\oC_{\q}}

\usepackage[curve,matrix,arrow]{xy} 

\usepackage[bookmarks=false, hyperfigures=false, a4paper]{hyperref}

\newcommand{\Dx}{\dd_x}
\newcommand{\Dy}{\dd_y}

\newcommand{\Dp}[1]{\,_{\phantom{h}}^{\underline{#1\!}\,}}
\newcommand{\dx}{\upxi}
\newcommand{\dy}{\upeta}

\newcommand{\Le}{\boldsymbol{\mathscr{E}}}
\newcommand{\Lh}{\boldsymbol{\mathscr{H}}}
\newcommand{\Lf}{\boldsymbol{\mathscr{F}}}

\newcommand{\dpplane}{\underline{\oC}{}_{\q}[x,y]} 
\newcommand{\dpshort}{\underline{\oC}{}_{\q}}
\newcommand{\dpOmega}{\underline{\Omega}_{\q}}

\newcommand{\rep}{\mathscr}
\newcommand{\repX}{\rep{X}}

\newcommand{\repP}{\rep{P}}

\newcommand{\repW}{\rep{W}}

\newcommand{\bref}[1]{\textbf{\ref{#1}}}

\newcommand{\SL}[1]{s\ell(#1)}

\renewcommand{\geq}{\geqslant}
\renewcommand{\leq}{\leqslant}

\newcommand{\tensor}{\otimes}

\newcommand{\q}{\mathfrak{q}}

\newcommand{\U}{\overline{\mathscr{U}}_{\q} s\ell(2)}
\newcommand{\UL}{\boldsymbol{{\mathscr{U}}}_{\q} s\ell(2)}

\newcommand{\D}{\mathscr{D}}

\newcommand{\ffrac}[2]{\raisebox{.5pt}{\mbox{\footnotesize$\displaystyle\frac{#1}{#2}$}}}
\newcommand{\fffrac}[2]{\raisebox{.9pt}{\mbox{\scriptsize$\displaystyle\frac{#1}{#2}$}}}

\newcommand{\qfac}[1]{[#1]!\,}
\newcommand{\qint}[1]{[#1]}

\newcommand{\qbin}[2]{\mathchoice%
  {\qbinm{#1}{#2}}{\qbinmm{#1}{#2}}%
  {\qbinmm{#1}{#2}}{\qbinmm{#1}{#2}}}
\newcommand{\qbinm}[2]{\mbox{\footnotesize$\displaystyle
    \genfrac{[}{]}{0pt}{}{#1}{#2}$}}
\newcommand{\qbinmm}[2]{\genfrac{[}{]}{0pt}{}{#1}{#2}}

\newcommand{\nbin}[2]{\mathchoice%
  {\nbinm{#1}{#2}}{\nbinmm{#1}{#2}}%
  {\nbinmm{#1}{#2}}{\nbinmm{#1}{#2}}}
\newcommand{\nbinm}[2]{\mbox{\footnotesize$\displaystyle
    \genfrac{(}{)}{0pt}{}{#1}{#2}$}}
\newcommand{\nbinmm}[2]{\genfrac{(}{)}{0pt}{}{#1}{#2}}

\newcommand{\dd}{\partial}

\newcommand{\oC}{\mathbb{C}}
\newcommand{\oP}{\mathbb{P}}
\newcommand{\oN}{\mathbb{N}}

\numberwithin{equation}{section}
\makeatletter
\@addtoreset{equation}{section}
\@addtoreset{subsubsection}{section}

\def\@secnumfont{\bfseries}
\def\subsubsection{\@startsection{subsubsection}{3}%
  \z@{.5\linespacing\@plus.7\linespacing}{-.5em}%
  {\normalfont\bfseries}}
\def\paragraph{\@startsection{paragraph}{4}%
  \z@\z@{-\fontdimen2\font}%
  \normalfont\bfseries}
\def\subparagraph{\@startsection{subparagraph}{5}%
  \z@\z@{-\fontdimen2\font}%
  \normalfont\bfseries}

\makeatother

\swapnumbers

\newtheorem{Lemma}[subsection]{Lemma}
\newtheorem{lemma}[subsubsection]{Lemma}

\theoremstyle{definition}

\begin{document}
\vspace*{-1.5\baselineskip}

\title[Divided-power quantum plane]{Quantum-$s\ell(2)$ action on a
  divided-power quantum plane at even roots of unity}

\author[Semikhatov]{A.M.~Semikhatov}%

\address{\mbox{}\kern-\parindent Lebedev Physics Institute
  \hfill\mbox{}\linebreak \texttt{ams@sci.lebedev.ru}}

\begin{abstract}
  We describe a nonstandard version of the quantum plane, the one in
  the basis of divided powers at an even root of unity
  $\q=e^{i\pi/p}$.  It can be regarded as an extension of the ``nearly
  commutative'' algebra $\oC[X,Y]$ with $X Y =(-1)^p Y X$ by
  nilpotents.  For this quantum plane, we construct a
  Wess--Zumino-type de~Rham complex and find its decomposition into
  representations of the $2p^3$-dimensional quantum group $\U$ and its
  Lusztig extension $\UL$; the quantum group action is also defined on
  the algebra of quantum differential operators on the quantum plane.
\end{abstract}

\maketitle

\thispagestyle{empty}

\section{Introduction}
Recent studies of logarithmic conformal field theory models have shown
the remarkable fact that a significant part of the structure of these
models is captured by factorizable ribbon quantum groups at roots of
unity~\cite{[FGST],[FGST2],[FGST3],[FGST-q]}.  This fits the general
context of the Kazhdan--Lusztig correspondence$/$duality between
conformal field theories (vertex-operator algebras) and quantum
groups~\cite{[KL]}, but in logarithmic conformal field theories this
correspondence goes further, to the extent that the relevant quantum
groups may be considered ``\textit{hidden symmetries}'' of the
corresponding logarithmic models (see~\cite{[S-q]} for a review
and~\cite{[GT],[S-U],[BFGT]} for further development).  This motivates
further studies of the quantum groups that are Kazhdan--Lusztig-dual
to logarithmic conformal field theories.  For the class of $(p,1)$
logarithmic models (with $p=2,3,\dots$), the dual quantum group is
$\U$ generated by $E$, $K$, and $F$ with the relations
\begin{gather}\label{the-qugr}
  \begin{aligned}
    KEK^{-1}&=\q^2E,\quad
    KFK^{-1}=\q^{-2}F,\\
    [E,F]&=\ffrac{K-K^{-1}}{\q-\q^{-1}},
  \end{aligned}
  \\
  \label{the-constraints}
  E^{p}=F^{p}=0,\quad K^{2p}=1
\end{gather}
at the even root of unity\enlargethispage{1\baselineskip}
\begin{gather}\label{the-q}
  \q=e^{\frac{i\pi}{p}},
\end{gather}
This $2p^3$-dimensional quantum group first appeared in~\cite{[AGL]}
(also see~\cite{[FHT],[Ar]}; a somewhat larger quantum group was
studied in~\cite{[Erd]}).

In this paper, we study not the $\U$ quantum group itself but an
algebra carrying its action, a version of the quantum
plane~\cite{[M-Montr]}.\footnote{Quantum planes and their relation to
  $s\ell(2)$ quantum groups is a ``classic'' subject discussed in too
  many papers to be mentioned here; in addition to~\cite{[M-Montr]},
  we just note~\cite{[WZ],[M]}, and~\cite{[LR93]}.  In a context close
  to the one in this paper, the quantum plane was studied at the third
  root of unity in~\cite{[Coq],[DNS],[CGT]}.}  Quantum
planes\pagebreak[3] (which can be defined in slightly different
versions at roots of unity, for example, infinite or finite) have the
nice property of being \textit{$\mathscr{H}$-module algebras} for
$\mathscr{H}$ given by an appropriate version of the quantum
$s\ell(2)$ (an $\mathscr{H}$-module algebra is an associative algebra
endowed with an $\mathscr{H}$ action that ``agrees'' with the algebra
multiplication in the sense that $X\,(u\,v)=\sum(X'u)\,(X''v)$ for all
$X\in\mathscr{H}$).  Studying such algebras is necessary for extending
the Kazhdan--Lusztig correspondence with logarithmic conformal field
theory, specifically, extending it to the level of fields; module
algebras are to provide a quantum-group counterpart of the algebra of
fields in logarithmic models when these are described in a manifestly
quantum-group invariant way.  Another nice feature of quantum planes,
in relation to the Kazhdan--Lusztig correspondence in particular, is
that they allow a ``covariant calculus''~\cite{[WZ]}, i.e., a de~Rham
complex with the space of $0$-forms being just the quantum plane and
with a differential that commutes with the quantum-group action.  (A
``covariant calculus'' on a differential $\U$-module algebra was also
considered in~\cite{[S-U]} in a setting ideologically similar to but
distinct from a quantum plane.)

Here, we explore a possibility that allows extending the $\U$-action
on a module algebra to the action of the corresponding Lusztig quantum
group $\UL$, with the additional $s\ell(2)$ generators conventionally
written as $\Le=\fffrac{E^p}{[p]!}$ and
$\Lf=\fffrac{F^p}{[p]!}$.\footnote{We use the standard notation
  $\qint{n}=\fffrac{\q^n - \q^{-n}}{\q - \q^{-1}}$,
  $\qfac{n}=\qint{1}\qint{2}\dots\qint{n}$, and
  $\mbox{\scriptsize$\displaystyle\genfrac{[}{]}{0pt}{}{m}{n}$}{}
  =\fffrac{\qfac{m}}{\qfac{m-n}\qfac{n}}$}$^{,}$\footnote{The
  counterpart of this $s\ell(2)$ Lie algebra on the ``logarithmic''
  side of the Kazhdan--Lusztig correspondence, in particular,
  underlies the name ``triplet'' for the extended chiral algebra of
  the $(p,1)$ logarithmic models, introduced
  in~\cite{[K-first],[GK1],[GK2]} (also see~\cite{[G-alg]}; the
  triplet algebra was shown to be defined in terms of the kernel of a
  screening for general~$p$ in~\cite{[FHST]} and was studied, in
  particular, in~\cite{[CF],[EF],[AM-triplet]}).} We use Lusztig's
trick of divided powers twice: to extend the $\U$ quantum group and to
``distort'' the standard quantum plane $\oCq[x,y]$, i.e., the
associative algebra on $x$ and $y$ with the relation
\begin{gather}\label{xyyx}
  y\,x = \q\,x\,y,
\end{gather}
by passing to the divided powers
\begin{equation*}
  x\Dp{m}=\ffrac{x^m}{[m]!},\quad
  y\Dp{m}=\ffrac{y^m}{[m]!},\qquad m\geq0.
\end{equation*}
Constructions with divided powers~\cite{[L-q1],[L]} are interesting at
roots of unity, of course; in our root-of-unity case~\eqref{the-q},
specifically, $\qint{p}=0$, and the above formulas cannot be viewed as
a change of basis.  We actually define \textit{the divided-power
  quantum plane} $\dpplane$ (or just $\dpshort$ for brevity) to be the
span of $x\Dp{n}\,y\Dp{m}$, $m,n\geq0$, with the (associative)
multiplication
\begin{equation}\label{xxdiv}
  x\Dp{m} x\Dp{n}=\qbin{m+n}{m}\,x\Dp{m+n},\quad
  y\Dp{m} y\Dp{n}=\qbin{m+n}{m}\,y\Dp{m+n}
\end{equation}
and relations
\begin{equation}
  y\Dp{n}\,x\Dp{m} = \q^{n m} x\Dp{m}\, y\Dp{n}.
\end{equation}
(Somewhat more formally, \eqref{xxdiv} can be considered relations as
well.)  Divided-power quantum spaces were first considered
in~\cite{[Hu]}, in fact, in a greater generality (in an arbitrary
number of dimensions, correspondingly endowed with an $\mathscr{U}_q
s\ell(n)$ action).\footnote{I thank N.~Hu for pointing
  Ref.~\cite{[Hu]} out to me and for useful remarks.}  Our quantum
plane carries a $\U$ action that makes it a $\U$-module algebra (a
particular case of the construction in~\cite{[Hu]}),
\begin{equation}\label{U-acts-div}
  \begin{split}
    E(x\Dp{m} y\Dp{n}) &=  \qint{m+1} x\Dp{m + 1} y\Dp{n - 1},\\
    K(x\Dp{m} y\Dp{n}) &= \q^{m-n} x\Dp{m} y\Dp{n},\\
    F(x\Dp{m} y\Dp{n}) &= \qint{n+1}\,x\Dp{m - 1} y\Dp{n + 1},
  \end{split}
\end{equation}
with $x\Dp{m}=0=y\Dp{m}$ for $m<0$.\footnote{In speaking of a
  $\U$-module algebra, we refer to the Hopf algebra structure on $\U$
  given by the comultiplication $\Delta(E)= E\tensor K + 1\tensor E$,
  $\Delta(K)=K\tensor K$, $\Delta(F)=F\tensor 1 + K^{-1}\tensor F$,
  counit $\epsilon(E)=\epsilon(F)=0$, $\epsilon(K)=1$, and antipode
  $S(E)=-E K^{-1}$, $S(K)=K^{-1}$, $S(F)=-K F$~\cite{[FGST]}.}  A
useful way to look at $\dpshort$ is to consider the ``almost
commutative'' polynomial ring $\oC_{\varepsilon}[X,Y]$ with $X
Y=\varepsilon\,Y X$ for $\varepsilon=(-1)^p$ and extend it with the
algebra of ``infinitesimals'' $x\Dp{m} y\Dp{n}$, $0\leq m,n\leq p-1$.
Then, in particular, $X$ and $Y$ behave under the $s\ell(2)$ algebra
of $\Le$ and $\Lf$ ``almost'' (modulo signs for odd~$p$) as
homogeneous coordinates on~$\oC\oP^1$.

We extend $\dpshort$ to a differential algebra, the algebra of
differential forms $\dpOmega$, which can be considered a
(Wess--Zumino-type) de~Rham complex of $\dpshort$, and then describe
the $\U$ action on $\dpOmega$ and find how it decomposes into
$\U$-representations (Sec.~\ref{sec:dp-plane}).  Further, we extend
the quantum group action on the quantum plane and differential forms
to a $\UL$ action.  In Sec.~\ref{sec:qdiff}, we introduce quantum
differential operators on $\dpshort$, make them into a module algebra,
and illustrate their use with a construction of the projective
quantum-group module~$\repP^+_1$.

\section{Divided-power quantum plane}\label{sec:dp-plane}
The standard quantum plane $\oC_q[x,y]$ at a root of unity and the
divided-power quantum plane $\dpshort$ at a root of unity can be
considered two different root-of-unity limits of the ``generic''
quantum plane.  The formulas below therefore follow from the
``standard'' ones (for example, $E(x^m y^n) =[n]x^{m+1} y^{n-1}$ and
$F(x^m y^n) =[m]x^{m-1} y^{n+1}$ for the quantum-$s\ell(2)$ action) by
passing to the divided powers at a generic $q$, when
$x\Dp{m}=\fffrac{x^m}{[m]!}$ \textit{is} a change of basis, and
setting $q$ equal to our $\q$ in~\eqref{the-q} in the end.  In
particular, we apply this simple strategy to the Wess--Zumino calculus
on the quantum plane~\cite{[WZ]}.

\subsection{De~Rham complex of $\dpplane$}
We define $\dpOmega$, the space of differential forms on $\dpshort$,
as the differential algebra on the $x\Dp{m}\,y\Dp{n}$ and $\dx$, $\dy$
with the relations (in addition to those in~$\dpshort$)
\begin{gather*}
  \dx\dy =-\q\dy\dx,\quad \dx^2=0,\quad\dy^2=0,
\end{gather*}
and
\begin{alignat*}{2}
  \dx\, x\Dp{n} &= \smash[t]{\q^{2 n}
    x\Dp{n} \dx},&\quad
  \dy\,x\Dp{n} &= \q^n x\Dp{n} \dy + \q^{2n-2}(\q^2 - 1)x\Dp{n - 1}
  y\dx,
  \\
  \dx\,y\Dp{n} &= \q^n y\Dp{n}\dx,&\quad
  \dy\, y\Dp{n} &= \q^{2 n} y\Dp{n} \dy
\end{alignat*}
and with a differential $d$ ($d^2=0$) acting as
\begin{align*}
  d (x\Dp{m}y\Dp{n}) &= \q^{m + n - 1}x\Dp{m - 1} y\Dp{n} \dx + \q^{n
    - 1}x\Dp{m} y\Dp{n - 1}\dy.
\end{align*}
In particular, $d x=\dx$ and $d y = \dy$.

The $\U$ action on $1$-forms that commutes with the differential is
given by
\begin{align*}
  E(x\Dp{m} y\Dp{n} \dx) &= \q [m+1] x\Dp{m + 1} y\Dp{n - 1} \dx,\\
  E(x\Dp{m} y\Dp{n} \dy) &= \q^{-1} [m+1] x\Dp{m + 1} y\Dp{n - 1} \dy
  + x\Dp{m} y\Dp{n} \dx,\\
  K(x\Dp{m} y\Dp{n} \dx) &=\q^{m-n+1}\,x\Dp{m} y\Dp{n} \dx,\\
  K(x\Dp{m} y\Dp{n} \dy) &=\q^{m-n-1}\,x\Dp{m} y\Dp{n} \dy,\\
  F(x\Dp{m} y\Dp{n} \dx) &= [n+1]\,x\Dp{m - 1} y\Dp{n + 1} \dx
  + \q^{n - m} x\Dp{m} y\Dp{n} \dy,\\
  F(x\Dp{m} y\Dp{n} \dy) &= [n+1]\,x\Dp{m - 1} y\Dp{n + 1} \dy
\end{align*}
(and on $2$-forms, simply by $E(x\Dp{m} y\Dp{n} \dx \dy) = [m+1]
x\Dp{m + 1} y\Dp{n - 1} \dx \dy$ and $F(x\Dp{m} y\Dp{n} \dx \dy) =
[n+1]x\Dp{m - 1} y\Dp{n + 1} \dx \dy$).

A feature of the divided-power quantum plane is that the monomials
$x\Dp{p m - 1}\,y\Dp{p n - 1}$ with $m,n\geq1$ are $\U$ invariants,
although not ``constants'' with respect to~$d$.  Also, the kernel of
$d$ in $\dpshort$ is given by $\oC 1$, and the cohomology of $d$
in $\dpOmega^1$ is zero.  A related, ``compensating,'' difference from
the standard quantum plane at the same root of unity amounts to zeros
occurring in the multiplication table: for $1\leq m,n\leq p-1$, it
follows that $x\Dp{m}\,x\Dp{n}=0$ whenever $m+n\geq p$).

\subsection{``Frobenius'' basis}\label{sec:Frobenius}
An ideal (the nilradical) is generated in $\dpshort$ from the
monomials $x\Dp{m}y\Dp{n}$ with $1\leq m,n\leq p-1$.  We now introduce
a new basis by splitting the powers of $x$ and $y$ accordingly
(actually, by passing to ``\textit{nondivided}'' powers of $x\Dp{p}$
and $y\Dp{p}$).

\subsubsection{}\label{sec:Frob-basis}
The divided-power quantum plane $\dpplane$ can be equivalently viewed
as the linear span of
\begin{equation}\label{Frob-basis}
  X^M\,Y^N\,x\Dp{m}\,y\Dp{n},\qquad
  0\leq m,n\leq p-1,\quad M,N\geq0,
\end{equation}
where $X=x\Dp{p}$ and $Y=y\Dp{p}$; the change of basis is explicitly
given by\footnote{Here and in what follows, we ``resolve'' the
  $q$-binomial coefficients using the general formula~\cite{[L]}
  \begin{gather*}
    \qbin{M p + m}{N p + n} = (-1)^{(M-1)N p + m N - n M}\,
    \nbin{M}{N}\qbin{m}{n}
  \end{gather*}
  for $0\leq m,n\leq p-1$ and $M\geq1$, $N\geq0$.}
\begin{equation*}
  x\Dp{M p + m}=(-1)^{M m}x\Dp{M p}\,x\Dp{m}=(-1)^{M m +
    \frac{M(M-1)}{2}p}\fffrac{1}{M!}\,X^M\,x\Dp{m}
\end{equation*}
and similarly for $y\Dp{N p + n}$ (cf.\ Proposition~2.4
in~\cite{[Hu]}).  Clearly, monomials in the right-hand side here are
multiplied ``componentwise,''
\begin{gather*}
  (X^M x\Dp{m})(X^N x\Dp{n})=X^{M+N}\qbin{m+n}{m} x\Dp{m+n}.
\end{gather*}

The relations involving $X$ and $Y$ are
\begin{alignat*}{2}
  X\,Y&=(-1)^p\,Y\,X,\kern-40pt\\
  X\,x\Dp{n}&=x\Dp{n}\,X,&\qquad X\,y\Dp{n}&=(-1)^n\,y\Dp{n}\,X,\\
  Y\,x\Dp{n}&=(-1)^n\,x\Dp{n}\,Y,&\qquad Y\,y\Dp{n}&=y\Dp{n}\,Y.
\end{alignat*}
We continue writing $\dpshort$ for $\oC_{\q}[X,Y,x,y]$ with all the
relations understood.  The quotient of $\dpshort$ by the nilradical is
the ``almost commutative'' polynomial ring $\oC_{\varepsilon}[X,Y]$,
where $X Y=\varepsilon\,Y X$ for $\varepsilon=(-1)^p$.

In $\dpOmega$, the relations involving $X$ and $Y$ are
\begin{alignat*}{2}
  \dx\,X&=X\,\dx,&\quad
  \dy\,X^M &= (-1)^M\,X^M\,\dy + (-1)^{M} M(\q^{-2} - 1)X^{M - 1}
  x\Dp{p-1}\,y\,\dx,
  \\
  \dx\,Y^{N} &= (-1)^N\, Y^{N}\,\dx,&
  \quad
  \dy\, Y&= Y\,\dy,
\end{alignat*}
and the differential is readily expressed as
\begin{multline*}
  d(X^M Y^N x\Dp{m} y\Dp{n}) =\q^{m+n-1}
  \begin{cases}
    X^M Y^N x\Dp{m-1} y\Dp{n} \dx,&m\neq0,\\
    -(-1)^{N p} M\, X^{M-1} Y^N x\Dp{p-1} y\Dp{n}\dx,& m=0
  \end{cases}\\*
  {}+\q^{n-1}
  \begin{cases}
    X^M Y^N x\Dp{m} y\Dp{n-1} \dy,& n\neq0,\\
    -(-1)^m N\, X^M Y^{N-1} x\Dp{m} y\Dp{p-1} \dy,& n=0.
  \end{cases}
\end{multline*}
In particular,
$d(X^M\,Y^N) = -(-1)^{N p}M\q^{-1} X^{M-1} Y^N x\Dp{p-1}\dx -N\q^{-1}
X^M Y^{N-1} y\Dp{p-1}\dy$.%

\subsubsection{}\label{U-XYloc}
The $\U$ action in the new basis~\eqref{Frob-basis} becomes
\begin{align*}
  E(X^M Y^N x\Dp{m} y\Dp{n})&=
  \begin{cases}
    [m+1] X^M Y^N x\Dp{m+1} y\Dp{n-1},& n\neq0,\\
    (-1)^m N [m+1] X^M Y^{N-1} x\Dp{m+1} y\Dp{p-1},& n=0,
  \end{cases}
  \\
  K(X^M Y^N x\Dp{m} y\Dp{n})&=
  (-1)^{M+N}\q^{m-n} X^M Y^N x\Dp{m} y\Dp{n},
  \\
  F(X^M Y^N x\Dp{m} y\Dp{n})&=
  \begin{cases}
    (-1)^{M+N} [n+1] X^M Y^N x\Dp{m-1} y\Dp{n+1},& m\neq 0,\\
    M(-1)^{M-1+N(p-1)} [n+1] X^{M-1} Y^N x\Dp{p-1} y\Dp{n+1},& m=0,
  \end{cases}
\end{align*}
and on $1$-forms, accordingly,
\begin{align*}
  E(X^M Y^N x\Dp{m} y\Dp{n} \dx)&=
  \begin{cases}
    \q\, [m+1] X^M Y^N x\Dp{m+1} y\Dp{n-1} \dx,& n\neq0,\\
    (-1)^m N \q\, [m+1] X^M Y^{N-1} x\Dp{m+1} y\Dp{p-1} \dx,& n=0,
  \end{cases}
  \\
  E(X^M Y^N x\Dp{m} y\Dp{n} \dy)&=
  \begin{cases}
    \q^{-1}[m+1] X^M Y^N x\Dp{m+1} y\Dp{n-1} \dy,& n\neq0,\\
    (-1)^m N \q^{-1}[m+1] X^M Y^{N-1} x\Dp{m+1} y\Dp{p-1} \dy,& n=0
  \end{cases}\\*
  &\quad{}+X^M Y^N x\Dp{m} y\Dp{n} \dx,
  \\
  F(X^M Y^N x\Dp{m} y\Dp{n} \dx)&=
  \begin{cases}
    (-1)^{M+N} [n+1] X^M Y^N x\Dp{m-1} y\Dp{n+1} \dx,& m\neq 0,\\
    M(-1)^{M-1+N(p-1)} [n+1] X^{M-1} Y^N x\Dp{p-1} y\Dp{n+1} \dx,& m=0
  \end{cases}\\*
  &\quad{}+ (-1)^{M+N}\q^{n-m} X^M Y^N x\Dp{m} y\Dp{n} \dy,
  \\
  F(X^M Y^N x\Dp{m} y\Dp{n} \dy)&=
  \begin{cases}
    (-1)^{M+N} [n+1] X^M Y^N x\Dp{m-1} y\Dp{n+1} \dy,& m\neq 0,\\
    M(-1)^{M-1+N(p-1)} [n+1] X^{M-1} Y^N x\Dp{p-1} y\Dp{n+1} \dy,&
    m=0.
  \end{cases}
\end{align*}

\subsection{Module decomposition} 
We next find how $\dpOmega$ decomposes as a $\U$-module.  The $\U$
representations encountered in what follows are the irreducible
representations $\repX^\pm_r$ ($1\leq r\leq p$), indecomposable
representations $\repW^\pm_r(n)$ ($1\leq r\leq p-1$, $n\geq2$), and
projective modules $\repP^\pm_{p-1}$.  All of these are described
in~\cite{[FGST2]}; we only briefly recall from~\cite{[FGST],[FGST2]}
that $\U$ has $2p$ irreducible representations $\repX^\pm_r$, $1\leq
r\leq p$, with $\dim\repX^\pm_r=r$.  With the respective projective
covers denoted by $\repP^\pm_r$, the ``plus'' representations are
distinguished from the ``minus'' ones by the fact that tensor products
$\repX^+_r\tensor\repX^+_s$ decompose into the $\repX^+_{r'}$ and
$\repP^+_{r'}$~\cite{[S-q]} (and $\repX^+_1$ is the trivial
representation).  We also note that $\dim\repP^\pm_r=2p$ if $1\leq
r\leq p-1$.

We first decompose the space of $0$-forms, the quantum plane
$\dpshort$ itself; clearly, the $\U$ action restricts to each graded
subspace $(\dpshort)_i$ spanned by $x\Dp{m}y\Dp{n}$ with $m+n=i$,
$i\geq0$.  To avoid an unnecessarily long list of formulas, we
explicitly write the decompositions for $(\dpshort)_{(i)}$ spanned by
$x\Dp{m}y\Dp{n}$ with $0\leq m+n \leq i$; the decomposition of each
graded component is in fact easy to
extract.\enlargethispage{.5\baselineskip}
\begin{lemma}
  $\dpshort$ decomposes into $\U$-representations as follows:
  \begin{align}\label{decomp(p-1)}
    (\dpshort)_{(p-1)} &= \bigoplus_{r=1}^{p}\repX^+_r
    \\
    \intertext{\textup{(}where each
      $\repX^+_r\in(\dpshort)_{r-1}$\textup{)},}
    \label{decomp(2p-1)}
    (\dpshort)_{(2p-1)} &=
    (\dpshort)_{(p-1)}\oplus\smash[t]{\bigoplus_{r=1}^{p-1}}\repW^{-}_r(2)
    \oplus 2\repX^-_{p}\\
    \intertext{\textup{(}where each
      $\repW^{-}_r(2)\in(\dpshort)_{p+r-1}$ and
      $\repX^-_{p}\in(\dpshort)_{2p-1}$\textup{)}, and, in general,}
    \label{dpdecomp-general}
    (\dpshort)_{(n p-1)}&= (\dpshort)_{((n-1)
      p-1)}\oplus\bigoplus_{r=1}^{p-1}\repW^{\pm}_r(n) \oplus
    n\repX^{\pm}_p
  \end{align}
  with the $+$ sign for odd $n$ and $-$ for even $n$.
\end{lemma}

The proof is by explicit construction and dimension counting; the
general pattern of the construction is already clear from the lower
grades.  For $1\leq r\leq p$, the irreducible representations
$\repX^+_r$ are realized as
\begin{equation*}
  \xymatrix{%
    x\Dp{r-1}\ar^F@/^/[r]
    &x\Dp{r-2}\,y\ar^(.5)E@/^/[l]
    \ar^F@/^/[r]
    &\ \dots\ \ar^{E}@/^/[l]\ar^{F}@/^/[r]
    &y\Dp{r-1}.\ar^{E}@/^/[l]
  }
\end{equation*}
(the $E$ and $F$ arrows represent the respective maps up to nonzero
factors).  The $\repX^+_r$ thus constructed exhaust the space
$(\dpshort)_{(p-1)}$, which gives~\eqref{decomp(p-1)}.  To pass to the
next range of grades, Eq.~\eqref{decomp(2p-1)}, it is easiest to
multiply the above monomials with $X$ or $Y$ and then note that (up to
sign factors) this operation commutes with the $\U$ action except ``at
the ends'' in half the cases, as is expressed by the identities
\begin{small}%
  \begin{align*}
    X E(X^M Y^N x\Dp{m} y\Dp{n}) - E(X^{M+1} Y^N x\Dp{m} y\Dp{n}) &=0,\\
    Y E(X^M Y^N x\Dp{m} y\Dp{n}) - (-1)^{M p} E(X^{M} Y^{N+1} x\Dp{m}
    y\Dp{n}) &=
    -\delta_{n,0}(-1)^{M p + m} [m+1] X^M Y^N x\Dp{m+1} y\Dp{p-1},\\
    X F(X^M Y^N x\Dp{m} y\Dp{n}) + F(X^{M+1} Y^N x\Dp{m} y\Dp{n}) &=
    \delta_{m,0}(-1)^{M + N(p-1)}[n+1] X^M Y^N x\Dp{p-1} y\Dp{n+1},\\
    Y F(X^M Y^N x\Dp{m} y\Dp{n}) + (-1)^{M p}F(X^{M} Y^{N+1} x\Dp{m}
    y\Dp{n}) &=0.
  \end{align*}%
\end{small}%
It then immediately follows that for each $r=1,\dots,p-1$, the states
\begin{small}
  \begin{gather}\label{W(2)-details}
    \xymatrix@=14pt{%
      X x\Dp{r-1}\ar^F@/^/[r]    
      &\ \dots\ \ar^{E}@/^/[l]\ar^{F}@/^/[r]
      &X y\Dp{r-1}\ar^{E}@/^/[l]
      \ar^F[dr]
      &&&&Y x\Dp{r-1}\ar_E[dl]
      \ar^F@/^/[r]    
      &\ \dots\ \ar^{E}@/^/[l]\ar^{F}@/^/[r]
      &Y y\Dp{r-1}\ar^{E}@/^/[l]
      \\
      &&&x\Dp{p-1}y\Dp{r}\ar^F@/^/[r]    
      &\ \dots\ \ar^{E}@/^/[l]\ar^{F}@/^/[r]
      &x\Dp{r}y\Dp{p-1}\ar^{E}@/^/[l]
    }
  \end{gather}%
\end{small}%
realize the representation~\cite{[FGST2]}
\begin{equation}\label{W(2)}
  \raisebox{-\baselineskip}{\mbox{$\displaystyle\repW^{-}_r(2)={}$}}
  \xymatrix@=14pt{%
    \repX^-_{r}\ar[dr]&&\repX^-_{r}\ar[dl]\\
    &\repX^+_{p-r}&
  }
\end{equation}
These $\repW^{-}_r(2)$ fill the space $(\dpshort)_{(2p-2)}$.
Next, in the subspace $(\dpshort)_{2p-1}$, this picture degenerates
into a sum of two $\repX^-_{p}$, spanned by $X x\Dp{p-1},\dots,X
y\Dp{p-1}$ and $Y x\Dp{p-1},\dots,Y y\Dp{p-1}$, with the result
in~\eqref{decomp(2p-1)}.

The longer ``snake'' modules follow similarly; for example, in the
$\repW^+_r(3)$ module
\begin{equation*}
  \xymatrix@=14pt{%
    \repX^+_{r}\ar[dr]&&\repX^+_{r}\ar[dl]\ar[dr]&&\repX^+_{r}\ar[dl]\\
    &\repX^-_{p-r}&&\repX^-_{p-r}
  }
\end{equation*}
the leftmost state in the left $\repX^+_{r}$ is given by $X^2
x\Dp{r-1}$ and the rightmost state in the right $\repX^+_{r}$ is $Y^2
y\Dp{r-1}$.  These $\repW^+_r(3)$ with $1\leq r\leq p-1$ exhaust
$(\dpshort)_{(3p-2)}$,
and for $r=p$, the snake degenerates into $3\repX^+_{p}$,
and so on, yielding~\eqref{dpdecomp-general} in general.

\medskip

Next, to describe the decomposition of the $1$-forms $\dpOmega^1$, we
use the grading by the subspaces
$(\dpOmega^1)_{i}=(\dpshort)_{i}\dx+(\dpshort)_{i}\dy$, $i\geq0$, and
set $(\dpOmega^1)_{(i)}=(\dpshort)_{(i)}\dx+(\dpshort)_{(i)}\dy$
accordingly.  For those representations $\repX\in\dpshort$ that are
unchanged under the action of~$d$, we write $d\repX$ for their
isomorphic images in $\dpOmega^1$, to distinguish the $d$-exact
representations in the decompositions that follow.  
\begin{lemma}
  The space of $1$-forms $\dpOmega^1$ decomposes into
  $\U$-representations as follows:
  \begin{align}\label{decompO(p-1)}
    (\dpOmega^1)_{(p-1)}&=\bigoplus_{r=2}^{p} d\repX^+_{r}\oplus
    \bigoplus_{r=1}^{p-2}\repX^+_{r}\oplus \repP^+_{p-1}
    \\
    \intertext{\textup{(}with $d\repX^+_r\in(\dpOmega^1)_{r-2}$,
      $\repX^+_r\in(\dpOmega^1)_{r}$, and
      $\repP^+_{p-1}\in(\dpOmega^1)_{p-1}$\textup{)},}
    \label{decompO(2p-1)}
    (\dpOmega^1)_{(2p-1)} &=(\dpOmega^1)_{(p-1)}
    \oplus\bigoplus_{r=2}^{p-1} d\repW^-_r(2) \oplus 2d\repX^-_p
    \oplus\repX^+_p \oplus\bigoplus_{r=1}^{p-2}\repW^-_r(2) \oplus
    2\repP^-_{p-1}
    \\
    \intertext{\textup{(}with $d\repW^-_r(2)\in(\dpOmega^1)_{p+r-2}$,
      $d\repX^-_p\in(\dpOmega^1)_{2p-2}$,
      $\repX^+_p\in(\dpOmega^1)_{p}$,
      $\repW^-_r(2)\in(\dpOmega^1)_{p+r}$, and
      $\repP^-_{p-1}\in(\dpOmega^1)_{2p-1}$\textup{)},}
    \label{decompO(3p-1)}
    (\dpOmega^1)_{(3p-1)} &=(\dpOmega^1)_{(2p-1)} \oplus
    \smash[t]{\bigoplus_{r=2}^{p-1}} d\repW^+_r(3) \oplus 3d\repX^+_p \oplus
    2\repX^-_p \oplus\smash[t]{\bigoplus_{r=1}^{p-2}} \repW^+_r(3)
    \oplus3\repP^{+}_{p-1}
    \\
    \intertext{\textup{(}with $d\repW^+_r(3)\in(\dpOmega^1)_{2p+r-2}$,
      $d\repX^+_p\in(\dpOmega^1)_{3p-2}$,
      $\repX^-_p\in(\dpOmega^1)_{2p}$,
      $\repW^+_r(3)\in(\dpOmega^1)_{2p+r}$, and
      $\repP^{+}_{p-1}\in(\dpOmega^1)_{3p-1}$\textup{)}, and, in general,}
    \label{decomp-last}
    (\dpOmega^1)_{(n p-1)} &=(\dpOmega^1)_{((n-1)p-1)}\oplus
    \smash[t]{\bigoplus_{r=2}^{p-1}} d\repW^{\pm}_r(n) \oplus
    n\,d\repX^{\pm}_p\\*
    &\notag \qquad\qquad\qquad{}\oplus (n-1)\repX^{\mp}_p
    \oplus\bigoplus_{r=1}^{p-2}\repW^{\pm}_r(n) \oplus
    n\repP^{\pm}_{p-1},
  \end{align}
  with the upper signs for odd $n$ and the lower signs for even~$n$.
\end{lemma}

The proof is again by construction.  To begin with
$(\dpOmega^1)_{(p-1)}$ in~\eqref{decompO(p-1)}, we first note that is
contains $d\repX^+_r\approx\repX^+_{r}$ spanned by\pagebreak[3] $d x\Dp{r-1}$,
$d(x\Dp{r-2}\,y)$, $\dots$, $d y\Dp{r-1}$ for each $r=2,\dots,p$
(clearly, the singlet $\repX^+_1$ spanned by~$1$ for $r=1$ is
annihilated by~$d$).  Next, for each $r=1,\dots,p-2$, there is an
$\repX^+_r$ spanned by the $1$-forms
\begin{equation*}
  x\Dp{r-1}y\dx-\q[r]x\Dp{r}\dy,\;\dots,\;
  [i+1]x\Dp{r-i-1}y\Dp{i+1}\dx - [r-i]\q^{i+1}x\Dp{r-i}y\Dp{i}\dy,\;
  \dots,\;
  [r]y\Dp{r}\dx - \q^rx y\Dp{r-1}\dy
\end{equation*}
(with just the singlet $y\dx-\q x\dy$ for $r=1$); for $r=p-1$,
however, these states are in the image of~$d$ and actually constitute
the socle of the projective module $\repP^+_{p-1}$ realized as
\begin{small}
  \begin{gather}\label{p-1-proj}
    \xymatrix@=16pt{%
      &x\Dp{p-1}\dy \ar[dl]_(.5)E
      \ar^{F}@/^/[r]&\qquad\dots\qquad\ar^{E}@/^/[l]\ar^{F}@/^/[r] &x
      y\Dp{p-2}\dy\ar^{E}@/^/[l] \ar[dr]^(.6)F
      \ar^E[];[ddl]+<26pt,26pt>
      \\
      x\Dp{p-1}\dx \ar[dr]^(.7)F&&&& y\Dp{p-1}\dy \ar[dl]_E
      \\
      & d(x\Dp{p-1}y) \ar^{F}@/^/[r] \ar_(.3)E[uur]-<15pt,25pt>;[]
      &\dots\ar^{E}@/^/[l]\ar^{F}@/^/[r]& d(x y\Dp{p-1})
      \ar^{E}@/^/[l] }
  \end{gather}%
\end{small}%
This gives~\eqref{decompO(p-1)}.

In the next grade, we have $(\dpOmega^1)_{p} = d\repW^-_2(2) \oplus
\repX^+_p$, where $\repX^+_p$ is spanned by
\begin{equation*}
  x\Dp{p-1}y\dx
  ,\
  \ldots
  ,\
  [i+1]x\Dp{p-i-1}y\Dp{i+1}\dx - [i]\q^{i+1}x\Dp{p-i}y\Dp{i}\dy
  ,\
  \ldots
  ,\
  x y\Dp{p-1}\dy
\end{equation*}
(with $i=0$ at the left end and $i=p-1$ at the right end).  Next,
\begin{equation*}
  (\dpOmega^1)_{(2p-2)}
  =(\dpOmega^1)_{(p-1)} \oplus\repX^+_p \oplus\bigoplus_{r=2}^{p-1}
  d\repW^-_r(2) \oplus 2d\repX^-_p
  \oplus\bigoplus_{r=1}^{p-2}\repW^-_r(2),
\end{equation*}
where the $\repW^-_r(2)$ modules with $r=1,\dots,p-2$ are spanned by
(with nonzero factors dropped)
\begin{small}
  \begin{equation}\label{OW2}
    \mbox{}\kern-10pt
    \xymatrix@C=0pt{%
      \mbox{\footnotesize$\begin{array}{l}
          X x\Dp{r-1}y\dx\\
          {}-\q[r]X x\Dp{r}\dy
        \end{array},$}
      \dots,
      \mbox{\footnotesize$\begin{array}{l}
          [r]X y\Dp{r}\dx\\
          {}-\q^r X x y\Dp{r-1}\dy
        \end{array}$}\ar^(.5)F[]+<70pt,-20pt>;[dr]+<-50pt,14pt>
      &&
      \mbox{\footnotesize$\begin{array}{l}
          Y x\Dp{r-1} y\dx\\
          -\q[r]Y x\Dp{r}\dy
        \end{array},$}
      \dots,
      \mbox{\footnotesize$\begin{array}{l}
          [r]Y y\Dp{r}\dx\\
          -\q Y x y\Dp{r-1}\dy
        \end{array}$}
      \ar_(.5)E[]+<-70pt,-20pt>;[dl]+<50pt,14pt>
      \\
      &x\Dp{p-1}y\Dp{r+1}\dx,\ \dots,\
      x\Dp{r+1}y\Dp{p-1}\dy
    }
  \end{equation}%
\end{small}%
In the next grade, the two $\repP^-_{p-1}$ modules occurring
in~\eqref{decompO(2p-1)} are given by
\begin{small}%
  \begin{gather}\label{two-proj}
    \xymatrix@C=0pt{%
      &X x\Dp{p-1}\dy, \dots,X x y\Dp{p-2}\dy
      \ar[]+<-40pt,-12pt>;[dl]_(.5)E \ar[]+<40pt,-12pt>;[dr]^(.6)F
      \\
      X x\Dp{p-1}\dx \ar[];[dr]+<-50pt,12pt>^(.7)F&& X y\Dp{p-1}\dy
      \ar[];[dl]+<50pt,12pt>_E
      \\
      & d(X x\Dp{p-1}y), \dots, d(X xy\Dp{p-1}) } \xymatrix@C=0pt{%
      &Y x\Dp{p-1}\dy, \dots,Y x y\Dp{p-2}\dy
      \ar[]+<-40pt,-12pt>;[dl]_(.5)E \ar[]+<40pt,-12pt>;[dr]^(.6)F
      \\
      Y x\Dp{p-1}\dx \ar[];[dr]+<-50pt,12pt>^(.7)F&& Y y\Dp{p-1}\dy
      \ar[];[dl]+<50pt,12pt>_E
      \\
      & d(Y x\Dp{p-1} y),\dots, d(Y x y\Dp{p-1}) }
  \end{gather}%
\end{small}%
(the left--right arrows are again omitted for compactness).  Next, we
have
\begin{equation*}
  (\dpOmega^1)_{(3p-2)}
  =(\dpOmega^1)_{(2p-1)}\oplus 2\repX^-_p\oplus \bigoplus_{r=2}^{p-1}
  d\repW^+_r(3) \oplus
  3d\repX^+_p\oplus\bigoplus_{r=1}^{p-2}\repW^+_r(3), 
\end{equation*}
where,
for example, each $\repW^+_r(3)$ is spanned by monomials that can be
easily constructed starting with $X^2 x\Dp{r-1}y\dx-\q[r]X^2
x\Dp{r}\dy$ in the top-left corner, and so on; the pattern readily
extends to higher degrees by multiplying with powers of $X$ and $Y$,
yielding~\eqref{decomp-last}.

\subsection{Extending to Lusztig's quantum group}
\subsubsection{}\label{sec:LGQ}
The above formulas for the $\U$ action of course imply that $E^p$ and
$F^p$ act on $\dpOmega$ by zero.  But using Lusztig's trick once
again, we can define the action of Lusztig's divided-power quantum
group $\UL$ that extends our $\U$.  As before, temporarily taking
the quantum group deformation parameter $q$ generic, evaluating the
action of $\Le=\fffrac{1}{[p]!}\,E^{p}$ and
$\Lf=\fffrac{1}{[p]!}\,F^{p}$, and finally setting $q=\q$ gives
\begin{align*}
  \Le(X^M\,Y^N\,x\Dp{m} y\Dp{n}) &=
  N(-1)^{(N-1)p + n + m}\,X^{M+1}\,Y^{N-1}\,x\Dp{m} y\Dp{n},\\
  \Lf(X^M\,Y^N\,x\Dp{m} y\Dp{n}) &=
  M(-1)^{(M-1)p}\,X^{M-1}\,Y^{N+1}\,x\Dp{m} y\Dp{n}.\\
  \intertext{Then, clearly, $\Lh=[\Le,\Lf]$ acts on $\dpshort$ as}
  \Lh(X^M\,Y^N\,x\Dp{m} y\Dp{n}) &=
  (-1)^{(M+N-1)p + m + n}(M - N) X^M\,Y^N\,x\Dp{m} y\Dp{n}.
\end{align*}
Modulo the ``slight noncommutativity'' of $X$ and $Y$ for odd $p$, we
here have the standard $\SL2$ action on functions of the homogeneous
coordinates on $\oC\oP^1$:
\begin{equation*}
  \Le f(X,Y,x,y) = f(X,Y,x,y)\,
  \overleftarrow{\!\ffrac{\dd}{\dd Y}\!}\,\,X,
  \quad
  \Lf f(X,Y,x,y) = Y\,\overrightarrow{\!\ffrac{\dd}{\dd X}\!}\,\,
  f(X,Y,x,y)
\end{equation*}
From the standpoint of this $\oC\oP^1$, the $x\Dp{m}$ and $y\Dp{n}$
are to be regarded as some kind of ``infinitesimals,'' almost (modulo
signs) invisible to $\mathscr{U}s\ell(2)$ generated by $\Le$ and
$\Lf$.

The $\Le$ and $\Lf$ act on $\dpOmega^1$ as
\begin{align*}
  \Le(X^M Y^N x\Dp{m} y\Dp{n} \dx) &=
  -N(-1)^{(N-1)p + n + m}\,X^{M+1}\,Y^{N-1}\,x\Dp{m} y\Dp{n}\,\dx,\\
  \Le(X^M Y^N x\Dp{m} y\Dp{n} \dy) &=
  -N(-1)^{(N-1)p + n + m}\,X^{M+1}\,Y^{N-1}\,x\Dp{m} y\Dp{n}\,\dy+{}\\*
  &\quad{}+\delta_{m,0}
  \begin{cases}
    N(-1)^n\,X^M Y^{N-1}\,x\Dp{p-1} y\Dp{n+1}\,\dx,& 0\leq n\leq p-2,\\
    X^M Y^N\,x\Dp{p-1}\,\dx,& n=p-1,
  \end{cases}
  \\
  \Lf(X^M Y^N x\Dp{m} y\Dp{n} \dx) &=
  M(-1)^{(M-1)p}\,X^{M-1}\,Y^{N+1}\,x\Dp{m} y\Dp{n}\,\dx\\*
  &{}+ \delta_{n,0}
  \begin{cases}
    M (-1)^{(M-1)p + m}\q^{-m-1} X^{M-1} Y^N x\Dp{m+1} y\Dp{p-1}\dy,&
    0\leq m\leq p-2,\\
    (-1)^{(M+N)p}\,X^M Y^N\,y\Dp{p-1}\,\dy,& m=p-1,
  \end{cases}
  \\
  \Lf(X^M Y^N x\Dp{m} y\Dp{n} \dy) &=
  M(-1)^{(M-1)p}\,X^{M-1}\,Y^{N+1}\,x\Dp{m} y\Dp{n}\,\dy.
\end{align*}

\subsubsection{$\UL$ representations}Using the above formulas, it is
easy to see how the $\U$ representations in
\eqref{decomp(p-1)}--\eqref{decomp-last} combine into
$\UL$-representations.  The simple result is that
decomposition~\eqref{dpdecomp-general} rewrites in terms of
$\UL$-representations as
\begin{gather*}
  (\dpshort)_{(n p-1)}= (\dpshort)_{((n-1)
    p-1)}\oplus\bigoplus_{r=1}^{p-1}\underbracket{\repW^{\pm}_r(n)}
  \oplus \underbracket{n\repX^{\pm}_p},
\end{gather*}
where the underbrackets indicate that each $\repW^{\pm}_r(n)$ becomes
an (indecomposable) $\UL$-representation and the $n$ copies of
$\repX^{\pm}_p$ are combined into an (irreducible)
$\UL$-repre\-sentation.\footnote{Recently, $\UL$-representations were
  systematically analyzed in~\cite{[BFGT]}; a more advanced notation
  than the underbracketing is used there.}  Similarly,
\eqref{decomp-last} becomes
\begin{gather*}
  (\dpOmega^1)_{(n p-1)} =(\dpOmega^1)_{((n-1)p-1)}\oplus
  \bigoplus_{r=2}^{p-1} \underbracket{d\repW^{\pm}_r(n)} \oplus
  \underbracket{n\,d\repX^{\pm}_p}\oplus \underbracket{(n-1)\repX^{\mp}_p}
  \oplus\bigoplus_{r=1}^{p-2}\underbracket{\repW^{\pm}_r(n)} \oplus
  \underbracket{n\repP^{\pm}_{p-1}},
\end{gather*}
where $(n-1)\repX^{\mp}_p$ is an irreducible, $\repW^{\pm}_r(n)$ an
indecomposable, and $n\repP^{\pm}_{p-1}$ a projective $\UL$-module.

To see this, we first note that $\Le$ and $\Lf$ act trivially
on~\eqref{decomp(p-1)}; each of the $\repX^+_r$ representations of
$\U$ is also an irreducible representation of $\UL$.  Next,
in~\eqref{decomp(2p-1)}, for the $\repW^{-}_r(2)$ modules shown
in~\eqref{W(2)}, $\Le$ and $\Lf$ map between the two $\repX^-_r$ and
(in this lowest, $n=2$, case of $\repW^{\pm}_r(n)$ modules) act
trivially on the socle,
\begin{equation}\label{W(2)-L}
  \xymatrix@=14pt{%
    \repX^-_{r}\ar[dr]
    \ar@{-->}@/^8pt/_{\Lf}[]+<4pt,10pt>;[rr]+<-4pt,10pt>
    &&\repX^-_{r}
    \ar@{-->}@/_16pt/_{\Le}[]+<0pt,12pt>;[ll]+<-0pt,12pt>
    \ar[dl]
    \\
    &\repX^+_{p-r}
    &
  }
\end{equation}
yielding an indecomposable representation of~$\UL$.  By the same
pattern, the sum of two $\repX^-_{p}$ in~\eqref{decomp(2p-1)} becomes
an irreducible $\UL$-representation.  This extends to the general
case~\eqref{dpdecomp-general}: $\Le$ and $\Lf$ act horizontally on all
the $\repW^{\pm}_r(n)$ modules and also on the $n$ copies of
$\repX^-_p$, $n\geq2$, making them into respectively indecomposable
and irreducible $\UL$-modules.

Passing to the space of $1$-forms, each of the $\U$-representations in
\eqref{decompO(p-1)} remains an $\UL$ representation; $\Le$ and $\Lf$
act trivially on the $\repX^+_{r}$ and interchange the ``corners'' of
$\repP^+_{p-1}$ shown in~\eqref{p-1-proj}:
\begin{equation*}
  \Lf(x\Dp{p-1}\dx)=y\Dp{p-1}\dy,\qquad
  \Le(y\Dp{p-1}\dy)=x\Dp{p-1}\dx
\end{equation*}
($\repP^+_{p-1}$ thus becomes a projective $\UL$-module).

Next, in~\eqref{decompO(2p-1)}, $\repX^+_p$ remains an irreducible
$\UL$-module (with $\Le$ and $\Lf$ acting trivially), and\pagebreak[3]
each of the $\repW^-_r(2)$ is an indecomposable $\UL$-module;
in~\eqref{OW2}, for example, $\Lf(X x\Dp{r-1}y\dx-\q[r]X x\Dp{r}\dy)=Y
x\Dp{r-1} y\dx -\q[r]Y x\Dp{r}\dy$, etc. Further, the sum of two
$\repP^-_{p-1}$ in~\eqref{decompO(2p-1)} becomes a projective
$\UL$-module; there, the top floor (see~\eqref{two-proj}) is a doublet
under the $s\ell(2)$ Lie algebra of $\Le$ and $\Lf$, the middle floor
is a singlet ($X y\Dp{p-1}\dy - (-1)^p Y x\Dp{p-1}\dx$) plus a triplet
($X x\Dp{p-1}\dx$, $X y\Dp{p-1}\dy +(-1)^p Y x\Dp{p-1}\dx$, $Y
y\Dp{p-1}\dy$),
and the lowest floor is also a doublet.

The organization of $\UL$ representations in the general case
in~\eqref{decomp-last} is now obvious; it largely follows from the
picture in the lower degrees by multiplying with powers of $X$ and
$Y$.

\subsubsection{``Lusztig--Frobenius localization''}
The formulas in~\bref{sec:Frobenius} can be extended from nonnegative
integer to integer powers of $X$ and $Y$.  This means passing from
$\dpshort=\oC_{\q}[X,Y,x,y]$ to $\oC_{\q}[X,X^{-1},Y,Y^{-1},x,y]$,
which immediately gives rise to infinite-dimensional $\UL$
representations, as we now indicate very briefly.

Taking the picture in~\eqref{W(2)-details} (with $1\leq r\leq p-1$)
and multiplying all monomials by $X^{-1}$ destroys the south-east $F$
arrow and the $\Lf$ arrows shown in~\eqref{W(2)-L}, yielding
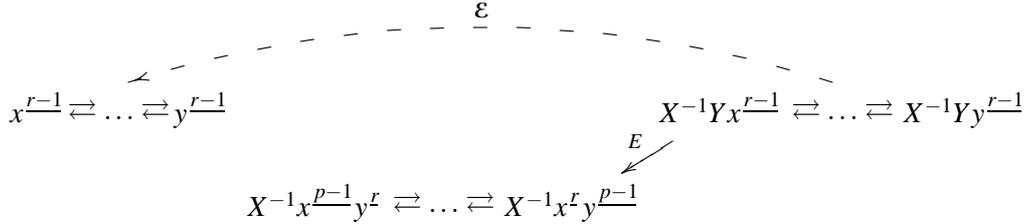
\begin{small}
  \begin{equation*}
    \xymatrix@=12pt@C=0pt{%
      x\Dp{r-1}
      &{\kern-6pt\rightleftarrows\ldots\rightleftarrows\kern-6pt}
      &y\Dp{r-1}
      &&&&X^{-1}Y x\Dp{r-1}\ar_E[dl]
      &*{\rightleftarrows\ldots\rightleftarrows}
      \ar@{-->}@/_20pt/_{\Le}[]+<-6pt,12pt>;[llllll]+<6pt,12pt>
      &X^{-1} Y y\Dp{r-1}
      \\
      &&&X^{-1} x\Dp{p-1}y\Dp{r}   
      &*{\rightleftarrows\ldots\rightleftarrows}
      &X^{-1} x\Dp{r}y\Dp{p-1}
    }
  \end{equation*}%
\end{small}%
But this now extends to the right, indefinitely, producing a pattern
shown in Fig.~\thefigure{} (the top diagram), where the $\Le$ and
$\Lf$ generators map between the respective monomials in each block.
A ``mirror-reflected'' picture follows by multiplying~\eqref{W(2)} by
$Y^{-1}$.  Multiplying by $X^{-1}Y^{-1}$ gives the bottom diagram in
Fig.~\thefigure{} (where the $\rightleftarrows$ arrows are omitted for
brevity).

After the extension by $X^{-1}$ and $Y^{-1}$, the differential $d$
acquires a nonzero cohomology in $\Omega^1\oC[X,X^{-1},Y,Y^{-1},x,y]$,
represented by $X^{-1}x\Dp{p-1}\dx$ and $Y^{-1}y\Dp{p-1}\dy$.

There is a curious possibility to further extend
$\oC[X,X^{-1},Y,Y^{-1},x,y]$ by adding a new element $\log(XY)$ that
commutes with $x$, $y$, $X$, and $Y$ and on which the differential and
the quantum group generators act, by definition, as follows:
\begin{alignat*}{2}
  d\log(XY) &=-\q^{-1} X^{-1} x\Dp{p-1}\dx -\q^{-1} Y^{-1}
  y\Dp{p-1}\dy,\kern-120pt
  \\
  E\log(XY) &= Y^{-1} x y\Dp{p-1},&\qquad \Le\log(XY) &= (-1)^p X
  Y^{-1},
  \\
  K\log(XY) &= \log(XY),\\
  F\log(XY) &= -X^{-1} x\Dp{p-1} y,& \Lf\log(XY) &= (-1)^p X^{-1} Y.
\end{alignat*}
\afterpage{%
  \rotatebox[origin=lB]{90}{\parbox{\textheight}{%
      \footnotesize\vspace*{-\baselineskip}
      \begin{gather*}
        \xymatrix@=12pt@C=4pt{%
          x\Dp{r-1}
          &*{\kern-6pt\rightleftarrows\ldots\rightleftarrows\kern-6pt}
          &y\Dp{r-1} &&&& X^{-1} Y x\Dp{r-1} \ar_E[dl]
          &*{\kern-6pt\rightleftarrows\ldots\rightleftarrows\kern-6pt}
          \ar@{-->}@/_20pt/_{\Le}[]+<-6pt,12pt>;[llllll]+<6pt,12pt>
          \ar@{-->}@/^20pt/^{\Lf}[]+<6pt,10pt>;[rrrrrr]+<0pt,10pt>
          &X^{-1} Y y\Dp{r-1} \ar[dr]^F &&&&X^{-2} Y^2 x\Dp{r}
          y\Dp{p-1}\ar[dl]^E &*{\kern-6pt\rightleftarrows\ldots}
          \ar@{-->}@/_16pt/^{\Le}[]+<0pt,6pt>;[llllll]+<6pt,6pt>
          \\
          &&&X^{-1}x\Dp{p-1}y\Dp{r}
          &*{\kern-6pt\rightleftarrows\ldots\rightleftarrows\kern-6pt}
          \ar@{-->}@/_20pt/_{\Lf}[]+<6pt,-12pt>;[rrrrrr]+<0pt,-12pt>
          &X^{-1}x\Dp{r}y\Dp{p-1} &&&&X^{-2} Y x\Dp{p-1} y\Dp{r}
          &*{\kern-6pt\rightleftarrows\ldots\rightleftarrows\kern-6pt}
          \ar@{-->}@/^16pt/_{\Le}[]+<0pt,-8pt>;[llllll]+<12pt,-8pt>
          &X^{-2} Y x\Dp{r} y\Dp{p-1} }
        \\[\baselineskip]
        \xymatrix@=12pt@C=2pt{%
          &&&Y^{-1} x\Dp{r-1} \ar@{-->}@/_/_{\Le}[lll]+<-20pt,0pt>
          \ar[dl]^E &*{\kern-6pt\ldots\kern-6pt} &Y^{-1} y\Dp{r-1}
          &&&&& X^{-1} x\Dp{r-1} &*{\kern-6pt\ldots\kern-6pt} &X^{-1}
          y\Dp{r-1} \ar@{-->}@/^/^{\Lf}[rrr]+<20pt,0pt> \ar[dr]^F &&&
          \\
          Y^{-2} x\Dp{p-1} y\Dp{r} &*{\kern-6pt\ldots\kern-6pt}
          &Y^{-2} x\Dp{r} y\Dp{p-1} &&&&X^{-1} Y^{-1} x\Dp{p-1}
          y\Dp{r} \ar@{-->}@/^20pt/_{\Le}[llllll]-<6pt,12pt>
          &*{\kern-6pt\ldots\kern-6pt} &X^{-1} Y^{-1} x\Dp{r}
          y\Dp{p-1}
          \ar@{-->}@/_20pt/^{\Lf}[]-<6pt,10pt>;[rrrrrr]+<16pt,-10pt>
          &&&&&X^{-2} x\Dp{p-1} y\Dp{r} &*{\kern-6pt\ldots\kern-6pt}
          &X^{-2} x\Dp{r} y\Dp{p-1} }
      \end{gather*}%
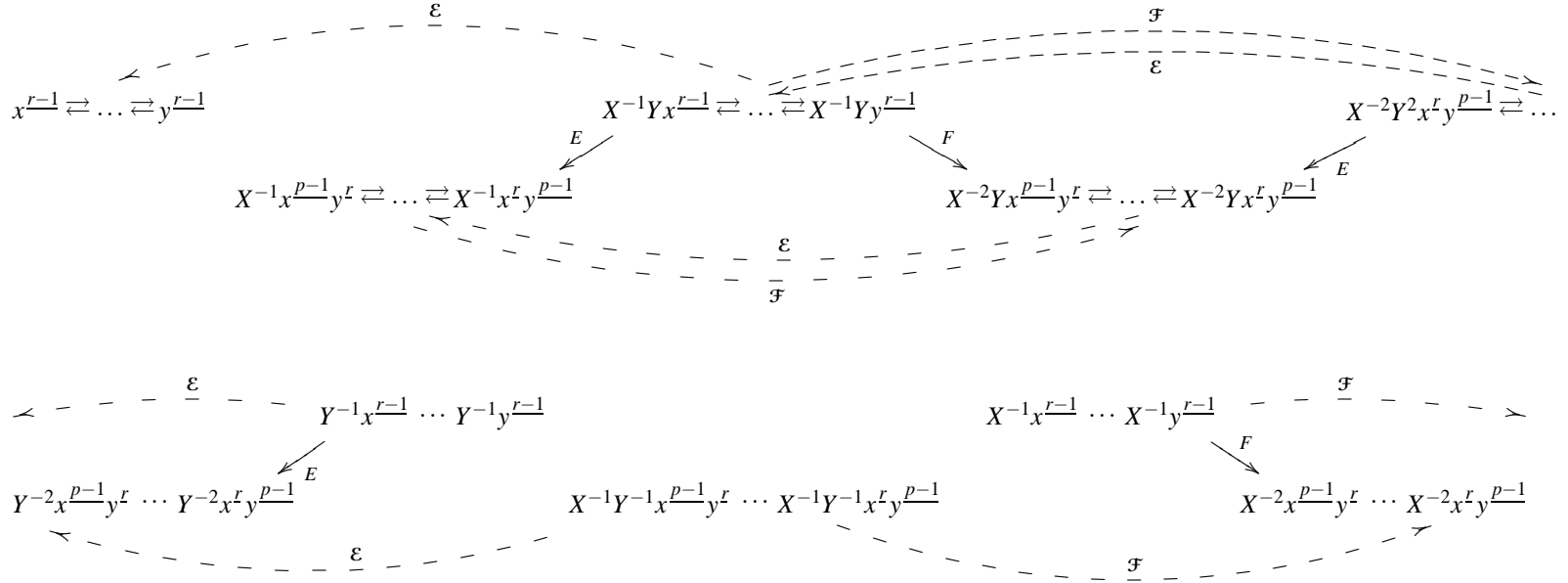
      \textsc{Figure}~\thefigure.  Two infinite-dimensional $\UL$
      representations realized on the quantum plane extended by
      negative powers of $X$ and $Y$.
      \\
      \addtocounter{figure}{1}%
      \begin{equation*}
        \mbox{}\kern20pt\xymatrix@=10pt@C=0pt{
          (-1)^p X Y^{-1}
          \ar^(.6){F}[dr]
          \ar^{d}[dd]
          &&&&\log(XY)
          \ar@{-->}@/_/_(.5){\Le}[llll]
          \ar@{-->}@/^/^(.5){\Lf}[rrrr]
          \ar_{E}[dl]
          \ar^{F}[dr]
          \ar^{d}[dd]
          &&
          &
          &(-1)^p X^{-1} Y
          \ar_(.6){E}[dl]
          \ar^{d}[dd]
          \\
          &-Y^{-1} x\Dp{p-1} y
          \ar^{d}[dd]
          &
          *{\kern-6pt\rightleftarrows\ldots\rightleftarrows\kern-6pt}
          &Y^{-1} x y\Dp{p-1}
          \ar^{d}[dd]
          &{}
          &-X^{-1} x\Dp{p-1} y
          \ar^{d}[dd]
          &
          *{\kern-2pt\rightleftarrows\ldots\rightleftarrows\kern-2pt}
          &(-1)^p X^{-1} x y\Dp{p-1}
          \ar^{d}[dd]
          \\
          \dots \ar^(.6){F}[dr] &&&& *{\begin{array}[t]{l}
              -\q^{-1} X^{-1} x\Dp{p-1}\dx\\[-3pt]
              {}-\q^{-1} Y^{-1} y\Dp{p-1}\dy
            \end{array}}
          \ar_{E}[dl]
          \ar^{F}[dr]
          &&&&\dots
          \ar^(.4){E}[];[dl]+<46pt,14pt>
          \\
          &
          *{\begin{array}[t]{l}
              \q^{-1} Y^{-1} x\Dp{p-2} y\dx\\[-3pt]
              {}-Y^{-1} x\Dp{p-1}\dy
            \end{array}}
          &
          *{\rightleftarrows\ldots\rightleftarrows}
          &
          *{\begin{array}[t]{l}
              -\q^{-1} Y^{-1} y\Dp{p-1}\dx\\[-3pt]
              {}-\q^{-2} Y^{-1} x y\Dp{p-2}\dy
            \end{array}}
          &{}
          &
          *{\begin{array}[t]{l}
              \q^{-1} X^{-1} x\Dp{p-2} y\dx\\[-3pt]
              {}-X^{-1} x\Dp{p-1}\dy
            \end{array}}
          &
          *{\rightleftarrows\ldots\rightleftarrows}
          &
          *{\begin{array}[t]{l}
              (-1)^{p-1}\q^{-1} X^{-1} y\Dp{p-1}\dx\\[-3pt]
              {}-(-1)^p \q^{-2} X^{-1} x y\Dp{p-2} \dy
            \end{array}}
          \\
          {}&
        }
      \end{equation*}%
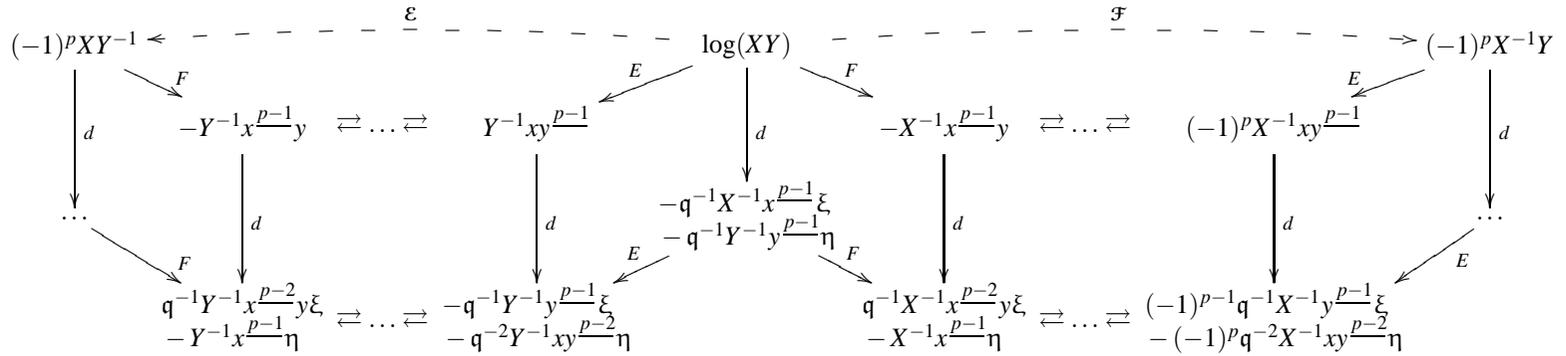
      \textsc{Figure}~\thefigure.  A $\UL$ representation involving
      $\log(XY)$ on the extended quantum plane, and its map by the
      differential.%
      \addtocounter{figure}{1}\mbox{}}}}%
We then have the diagram shown in
\addtocounter{figure}{1}Fig.~\thefigure.  \addtocounter{figure}{-1}%
There,\pagebreak[3] the $\rightleftarrows\ldots\rightleftarrows$
arrows represent maps modulo nonzero factors; for all the other
arrows, the precise numerical factors are indicated.  The pattern
extends infinitely both left and right.  In the central part, two more
maps that did not fit the picture are
\begin{gather*}
  \displaystyle\xymatrix@=12pt@C=80pt{ X Y^{-1}
    \ar@{-->}@/^/_{\Lf}[dr] & &X^{-1} Y \ar@{-->}@/_/^{\Le}[dl]
    \\
    &1 }
\end{gather*}

\section{Quantum differential operators}\label{sec:qdiff}
We next consider differential operators on the divided-power quantum
plane $\dpshort$.

\subsection{}
We define $\Dx$ and $\Dy:\dpshort\to\dpshort$ standardly, in
accordance with
\begin{equation*}
  d f= \dx \Dx f + \dy \Dy f
\end{equation*}
for any function $f$ of divided powers of $x$ and $y$.  It then
follows that\footnote{The rule to commute $\Dx$ and $\Dy$ inherited
  from the Wess--Zumino-type complex is not the one in~\cite{[Hu]}, as
  was also noted in that paper.}
\begin{align*}
  \Dx \Dy &= \q \Dy \Dx
  \\[-8pt]
  \intertext{and}
  \smash[t]{\Dx(x\Dp{m} y\Dp{n})} &=
  \smash[t]{\q^{-m - 2n + 1} x\Dp{m-1} y\Dp{n}},\\
  \Dy(x\Dp{m} y\Dp{n}) &= \q^{-m - n + 1} x\Dp{m} y\Dp{n-1}
  \\[-1.3\baselineskip]
\end{align*}
for $m,n\in\oN$; in terms of the basis in~\bref{sec:Frob-basis},
therefore,
\begin{align*}
  \Dx(X^M Y^N x\Dp{m} y\Dp{n}) &=
  \begin{cases}
    (-1)^N \q^{-m - 2n + 1} X^M Y^N x\Dp{m-1} y\Dp{n},& m\neq 0,\\
    M (-1)^{N(p-1) + 1} \q^{-2n+1} X^{M-1} Y^N x\Dp{p-1} y\Dp{n},&
    m=0,
  \end{cases}
  \\
  \Dy(X^M Y^N x\Dp{m} y\Dp{n}) &=
  \begin{cases}
    (-1)^{M}\q^{-m-n+1} X^M Y^N x\Dp{m} y\Dp{n-1},& n\neq0,\\
    N(-1)^M (-1)^{m-1}\q^{-m+1} X^M Y^{N-1} x\Dp{m} y\Dp{p-1},& n=0
  \end{cases}
\end{align*}
for $0\leq m,n\leq p-1$ and $M,N\geq 0$.

\subsection{}\label{sec:all-D-comm}
Let $\D$ denote the linear span of $x\Dp{m} y\Dp{n} \Dx^a \Dy^b$ with
$m,n,a,b\geq 0$.  The following commutation relations are easily
verified
($\ell\in\oN$):
\begin{equation}
  \label{eq:commute-linear}
  \begin{aligned}
    \Dx^\ell\,x &= \q^{-2 \ell}x\Dx^\ell +
    \q^{-\ell+1}\qint{\ell}\Dx^{\ell-1}
    -\q^{-2\ell+1}\qint{\ell}(\q-\q^{-1})y\Dx^{\ell-1}\Dy,\\
    \Dy^\ell\,y &= \q^{-2\ell}y \Dy^{\ell} + \q^{ -\ell + 1}
    \qint{\ell}\Dy^{\ell - 1}.
  \end{aligned}
\end{equation}
(Here and in similar relations below, $x$ is to be understood as the
operator of multiplication by $x$, etc.)  It then follows that
$\Dx^{p}$ and $\Dy^{p}$ are somewhat special: they (anti)commute with
all the divided powers $x\Dp{n}$ and $y\Dp{n}$ with $n<p$.  In
general, we have
\begin{multline*}
  \Dx^{\ell} x\Dp{m} y\Dp{n} =
  \q^{-2\ell m - \ell n}x\Dp{m} y\Dp{n} \Dx^{\ell}\\
  \shoveright{{}+ \sum_{i=1}^{\ell}\!\sum_{j=0}^{i} \q^{-2 m \ell -
      \ell n + i(\ell + m - n) - \frac{i(i - 1)}{2} - j(\ell + 1)}
    \qbin{\ell}{i}\qbin{i}{j}\qbin{n + j}{j} \qfac{j}\left(1 -
      \q^2\right)^j x\Dp{m - i} y\Dp{j + n} \Dx^{\ell - i} \Dy^{j},}
  \\
  \shoveleft{\Dy^{\ell} x\Dp{m} y\Dp{n} = \q^{-2\ell n - \ell
      m}x\Dp{m} y\Dp{n} \Dy^{\ell} + \sum_{i=1}^{\ell}\q^{-2 n \ell -
      \ell m + i(n + \ell) - \frac{i(i - 1)}{2}}
    \qbin{\ell}{i}x\Dp{m} y\Dp{n - i} \Dy^{\ell - i}}\\[-\baselineskip]
\end{multline*}
for all nonnegative integers $m$, $n$, and $\ell$, and hence the
commutation relations of $\Dx^p$ and $\Dy^p$ with elements of
$\dpshort$ can be written as
\begin{align*}
  \Dx^{p}\, X^M Y^N x\Dp{m} y\Dp{n} &= (-1)^{N p + n} X^M Y^N x\Dp{m}
  y\Dp{n} \Dx^{p} + M \q^{-\frac{p(p - 1)}{2}} X^{M-1} Y^N x\Dp{m}
  y\Dp{n},
  \\
  \Dy^{p}\,X^M Y^N x\Dp{m} y\Dp{n} &= (-1)^{Mp + m} X^M Y^N x\Dp{m}
  y\Dp{n} \Dy^{p} + (-1)^{M p} N \q^{-\frac{p(p - 1)}{2}} X^M Y^{N-1}
  x\Dp{m} y\Dp{n}
\end{align*}
for $M,N\geq1$ and $0\leq m,n\leq p-1$; in other words, $\Dx^p$ and
$\Dy^p$ can be represented as the (left) derivatives
\begin{equation*}
  \Dx^p=i^{1 - p}\ffrac{\dd}{\dd X},\qquad
  \Dy^p=i^{1 - p}\ffrac{\dd}{\dd Y}.
\end{equation*}
(Evidently, $\Dx^{p}\Dy^{p}=(-1)^{p}\Dy^{p}\Dx^{p}$.)  We here used
that $\q^{\frac{p(p - 1)}{2}} = i^{p - 1}$.

\begin{Lemma}
  The relations in $\D$ are compatible with the $\U$-module algebra
  structure if we define the $\U$ action as
  \begin{alignat*}{2}
    E \Dx &=-\q \Dy,&\qquad E \Dy &=0,\\
    K \Dx &= \q^{-1}\Dx,& K \Dy &= \q\Dy,\\
    F \Dx &=0, &\quad F \Dy &=-\q^{-1}\Dx.
  \end{alignat*}
\end{Lemma}
\noindent
The proof amounts to verifying that the relations
in~\bref{sec:all-D-comm} are mapped into one another under the $\U$
action.  Then $\Le$ and $\Lf$ act on the differential operators as
\begin{align*}
  \Le(\Dx^{Mp + m} \Dy^{N p + n}) &= -(-1)^{(N+1) p + n} M \Dx^{(M-1)
    p + m} \Dy^{(N+1) p + n},
  \\
  \Lf(\Dx^{M p + m} \Dy^{N p + n}) &= -(-1)^{(M+1) p + m} N \Dx^{(M+1)
    p + m}\Dy^{(N-1) p + n}.
\end{align*}
Modulo sign factors, this is the $s\ell(2)$ action on
$\oC[\Dx^p,\Dy^p]=\oC[\fffrac{\dd}{\dd X},\fffrac{\dd}{\dd Y}]$.

\subsection{A projective $\U$ module in terms of quantum differential
  operators} As an application of the $q$-differential operators on
$\dpshort$, we construct the projective $\U$-module $\repP^+_1$
(see~\cite{[FGST]}) as follows:
\begin{equation*}
  \xymatrix@C=4pt@R=12pt{
    &&&\pT
    \ar^{F}[dr]
    \ar_{E}[dl]
    \\  
    \pL_{p-1}
    &
    \kern-6pt\rightleftarrows \ldots \rightleftarrows\kern-6pt
    &\pL_1
    \ar_{F}[dr]
    &{}
    &\pR_1
    \ar^{E}[dl]
    &\kern-6pt\rightleftarrows \ldots \rightleftarrows\kern-6pt
    &
    \pR_{p-1}
    \\
    &&&
    \pB
  }
\end{equation*}
where the horizontal $\rightleftarrows$ arrows, as before, stand for
the action by $F$ and $E$ up to nonzero factors, and the actual
expressions for the module elements are as follows: first,\pagebreak[3]
\begin{gather*}
  \pT = -\q^{\frac{p(p - 1)}{2} + 1} \sum_{i=1}^p a_{i}\,x\Dp{p - i}
  y\Dp{i - 1} \Dx^{p - i} \Dy^{i -
    1},\quad\text{where}\quad
  a_i = \alpha + \smash[t]{\sum_{j=2}^i \ffrac{\q^{1-j}}{\qint{j-1}}},
\end{gather*}
then $\pL_i = E^i\pT$ and $\pR_i=F^i\pT$, with
\begin{alignat*}{2}
  \pL_{p-1} & = (-1)^{p-1} \qfac{p - 1} \q^{\frac{p(p - 1)}{2} + 1}
  x\Dp{p - 1}\Dy^{p - 1},&\ \pL_1 &= \q^{\frac{p(p - 1)}{2}}
  \!\sum_{i=1}^{p - 1}\!\q^{i + 2} x\Dp{p - i} y\Dp{i - 1}  \Dx^{p - i - 1} \Dy^{i},\\
  \pR_1 &= -\q^{\frac{p(p - 1)}{2}} \!\sum_{i=1}^{p - 1}\!\q^{i + 1}
  x\Dp{i - 1} y\Dp{p - i} \Dx^{i} \Dy^{p - 1 - i},&\ \pR_{p - 1}
  &= -\qfac{p - 1}\q^{\frac{p(p - 1)}{2} + 2} y\Dp{p - 1}\Dx^{p - 1}
\end{alignat*}
in particular, and, finally, $\pB = F\pL_1 = E\pR_1$ is
\begin{align*}
  \pB &=-\q^{\frac{p(p - 1)}{2} + 1} \sum_{i=1}^{p}x\Dp{p - i} y\Dp{i
    - 1} \Dx^{p - i} \Dy^{i - 1}.
\end{align*}
In the expression for $\pT$, $\alpha$ is an arbitrary constant
(clearly, adding a constant to all the $a_i$ amounts to redefining
$\pT$ by adding $\pB$ times this constant).  The normalization is here
chosen such that $\pB$ be a projector,
\begin{equation*}
  \pB  \pB = \pB.
\end{equation*}
(We note the useful identity $\qfac{p - 1}(\q - \q^{-1})^{p - 1}
=p\,\q^{\frac{p(p - 1)}{2}}$.)

The ``wings'' of the projective module 
are commutative,
\begin{equation*}
  \pL_i \pL_j=\pL_j \pL_i,\qquad
  \pR_i \pR_j=\pR_j \pR_i
\end{equation*}
for all $0\leq i,j\leq p-1$, where $\pL_0=\pR_0=\pB$, and, moreover,
$\pL_i \pL_j = \pR_i \pR_j=0$ whenever $i+j\geq p$.

A similar (but notably simpler) realization of $\repP^+_1$ in a
$\U$-module algebra of quantum differential operators on a ``quantum
line'' was given in~\cite{[S-U]}.

\section{Conclusion} 
Quantum planes provide a natural example of module algebras over
$s\ell(2)$ quantum groups (they do not allow realizing all of the
quantum-$s\ell(2)$ representations, but the corresponding quantum
differential operators make up a module algebra containing the
projective modules in particular).  By~\cite{[WZ]}, moreover,
$GL_q(2)$ can be \textit{characterized} as the ``quantum automorphism
group'' of the de~Rham complex of the quantum plane.  This is
conducive to the occurrence of quantum planes in various situations
where the $s\ell(2)$ quantum groups play a role.  From the standpoint
of the Kazhdan--Lusztig correspondence, the old subject of a quantum
$s\ell(2)$ action on the quantum plane is interesting in the case of
even roots of unity, which we detailed in this paper.

\subsubsection*{Acknowledgments}This paper was supported in part by
the RFBR grant 08-01-00737 and the grant LSS-1615.2008.2.


\parindent0pt

\end{document}